\documentclass[12pt,a4paper,reqno]{amsart}
\usepackage[margin=2.5cm]{geometry}
\usepackage[T1]{fontenc}
\usepackage{lmodern}
\usepackage[expansion=false]{microtype}
\usepackage{amsmath,amssymb,mathtools}
\usepackage{hyperref}
\usepackage{orcidlink}
\hypersetup{hidelinks,
 pdftitle={Exponential singularity and rank bounds for random symmetric matrices}}

\newtheoremstyle{paper}{6pt}{6pt}{\normalfont}{}{\bfseries}{.}{.5em}{}
\theoremstyle{paper}
\newtheorem{theorem}{Theorem}[section]
\newtheorem{proposition}[theorem]{Proposition}
\newtheorem{lemma}[theorem]{Lemma}
\newtheorem{corollary}[theorem]{Corollary}
\theoremstyle{remark}

\newtheorem*{remark*}{Remark}
\numberwithin{equation}{section}
\DeclareMathOperator{\rank}{rank}
\DeclareMathOperator{\corank}{corank}
\newcommand{\R}{\mathbb R}
\newcommand{\Prob}{\mathbb P}
\newcommand{\E}{\mathbb E}

\title{Exponential singularity and rank bounds for random symmetric matrices}
\author{Huaijin Liang}
\address{Institute for Financial Studies, Shandong University,
Jinan 250100, China}
\email{lianghuaij@gmail.com}
\author{Kexin Yu}
\address{School of Mathematics, Shandong University,
Jinan 250100, China}
\email{kexinyu01@gmail.com}
\author[Tingzhou Yu]{Tingzhou Yu\,\orcidlink{0009-0005-0414-7522}}
\address{Department of Mathematical and Statistical Sciences,
University of Alberta, Edmonton, AB T6G 2R3, Canada}
\email{tingzho1@ualberta.ca}


\date{}
\begin{document}
\raggedbottom
\begin{abstract}
Let $A$ be an $n\times n$ real symmetric random matrix whose entries above the diagonal are independent, each with largest atom at most $\beta\in (0,1)$.
The diagonal entries have an arbitrary joint distribution and are independent of the entries above the diagonal.
We prove that
\[
 \mathbb{P}[\mathrm{rank}(A)\le n-k]\le 2e^{-c_\beta nk},\qquad 1\le k\le n,
\]
where $c_\beta>0$ depends only on $\beta$.
This extends Han's corank bound to the full range $1\le k\le n$.
No moment assumptions are needed, and the entries need not be
centered or identically distributed.
In particular, the case $k=1$ gives an exponential singularity
bound under the uniform atom assumption.
Our argument also gives a new proof of the exponential singularity
bound for symmetric Rademacher matrices obtained by Campos,
Jenssen, Michelen, and Sahasrabudhe.
For every fixed $p\in (0,1)$, it yields the exponential singularity
bound conjectured by Costello and Vu for the adjacency matrix
of the Erd\H{o}s--R\'enyi random graph $G(n,p)$, where the case $p=1/2$ was previously proved by Han.
\end{abstract}

\maketitle

\section{Introduction}

We study the rank over $\R$ of $n\times n$ real symmetric random
matrices with independent entries above the diagonal.
Our main result requires only a uniform bound below one on the
largest atom of each of these entries and applies to all possible
coranks.

\begin{theorem}\label{thm:main}
Fix $\beta\in(0,1)$. There is a constant $c_\beta>0$, depending
only on $\beta$, such that the following holds for every positive
integer $n$.
Let $A=(a_{ij})_{1\le i,j\le n}$ be a real symmetric random matrix
whose entries above the diagonal are independent and satisfy
\[
 \sup_{x\in\R}\Prob[a_{ij}=x]\le\beta,
 \qquad 1\le i<j\le n.
\]
The diagonal entries have an arbitrary joint distribution.
The families $(a_{ii})_{i=1}^n$ and $(a_{ij})_{1\le i<j\le n}$
are independent.
Then, for every integer $1\le k\le n$,
\[
 \Prob[\rank(A)\le n-k]\le2e^{-c_\beta nk}.
\]
\end{theorem}

For matrices whose entries on and above the diagonal are independent
copies of a centered subgaussian random variable of variance one,
a bound of the form $e^{-cnk}$ for the probability of corank at least
$k$ was proved by Han in~\cite[Theorem~1.1]{Han} for
$1\le k\le c'\sqrt n$, where $c,c'>0$ depend only on the distribution.

Theorem~\ref{thm:main} removes the restriction on $k$.
It requires no moment assumptions, and the entries need not be
centered or identically distributed.
It also answers the question 
in~\cite[Remark~5.2]{Polavarapu}.
In particular,
\[
 \Prob[\det(A)=0]\le2e^{-c_\beta n}.
\]
Thus, an exponential singularity bound holds under the atom assumption. For symmetric Rademacher matrices, we provide a new proof of the bound of Campos,
Jenssen, Michelen, and Sahasrabudhe in \cite[Theorem 1.1]{CJMS}, see Section~\ref{sec:proofideas}.

Under a uniform atom bound and without moment assumptions,
polynomial singularity bounds for symmetric matrices whose entries above the diagonal may have different distributions were proved by Costello, Tao, and Vu
in~\cite[Theorem~6.5]{CTV} and by Manrique, P\'erez-Abreu, and Roy
in~\cite[Theorem~1(b) and Remark~2(b)]{MPR}.

For $\beta=1/2$, take the diagonal entries to be zero and the entries
above the diagonal to be independent Bernoulli$(1/2)$ variables.
The first $k$ rows are zero precisely when $k(n-k)+\binom{k}{2}$
independent entries are zero. For every $1\le k\le n$, this event has probability
\[
 2^{-k(n-k)-\binom{k}{2}}
 =2^{-kn+\binom{k+1}{2}}\ge2^{-nk}.
\]
This event forces $\rank(A)\le n-k$, so the exponent $nk$ is
optimal up to a constant.

For $p\in(0,1)$, let $G(n,p)$ be the random graph on $n$ vertices
in which each possible edge between distinct vertices is present
independently with probability $p$.

\begin{corollary}\label{cor:erdos-renyi}
Fix $p\in(0,1)$. There is a constant $c_p>0$, depending only on
$p$, such that the adjacency matrix
$A_G$ of $G(n,p)$ satisfies
\[
 \Prob[\rank(A_G)\le n-k]\le2e^{-c_pnk}
\]
for every positive integer $n$ and every integer $1\le k\le n$.
\end{corollary}

\begin{proof}
The diagonal of $A_G$ is zero, and its entries above the
diagonal are independent with largest atom $\max\{p,1-p\}<1$.
Apply Theorem~\ref{thm:main} with this atom bound.
\end{proof}

For the adjacency matrix of $G(n,p)$, a bound of this form was
proved by Han in~\cite[Theorem~1.3]{Han} for
$C_p\le k\le c\sqrt n$, where $C_p, c>0$ depend only on $p$;
when $p=1/2$, one may take $C_p=1$ (see also \cite[Remark 6.2]{Han}).
Corollary~\ref{cor:erdos-renyi} includes $k=1$ for every fixed
$p\in(0,1)$, giving $\Prob[\det(A_G)=0]\le2e^{-c_pn}$.
This establishes the conjectured exponential bound of Costello
and Vu in~\cite[Section~10]{CostelloVu} for every such $p$.
Theorem~\ref{thm:main} also applies to graphs with independent edges
whose probabilities are uniformly bounded away from zero and one,
with the constant depending only on these uniform bounds.

\begin{remark*}
The factor $2$ cannot be omitted if the bound is to hold in every
dimension. For example, take $\beta=1/2$ and
\[
 A=\begin{pmatrix}1&\xi\\\xi&1\end{pmatrix},
 \qquad \Prob[\xi=1]=\Prob[\xi=-1]=\frac12.
\]
Then $\rank(A)=1$ almost surely, so
$\Prob[\rank(A)\le1]=1>e^{-2c}$ for every $c>0$.
\end{remark*}

\subsection{Background}

The study of the singularity probability of random matrices goes
back to the work of Koml\'os in~\cite{Komlos}. For matrices with
independent Rademacher entries, an exponential upper bound was
obtained by Kahn, Koml\'os, and Szemer\'edi in~\cite{KKS}, and the sharp exponential rate was obtained by Tikhomirov \cite{Tikhomirov}. For symmetric matrices, Costello, Tao, and Vu \cite{CTV} proved that the singularity tends to zero, and successive improvements were obtained by Nguyen \cite{NguyenSymmetric}, Vershynin \cite{Vershynin}, Ferber and Jain \cite{FerberJain}, Campos,
Mattos, Morris, and Morrison \cite{CMMM}, and Campos,
Jenssen, Michelen, and Sahasrabudhe in \cite{CamposRevisited}.
The exponential bound for symmetric Rademacher matrices was
proved by Campos, Jenssen, Michelen, and Sahasrabudhe
in~\cite[Theorem~1.1]{CJMS}.
An exponential singularity bound was also proved
in~\cite[Theorem~1.1]{CJMSLeast}
for symmetric matrices whose entries on and above the diagonal
are independent copies of a centered subgaussian random variable
of variance one.

For a symmetric matrix $A$ whose entries on and above the diagonal
are independent and centered, each with variance one and subgaussian
norm at most $T\ge1$, the bound
$\Prob[\rank(A)\le n-k]\le e^{-cnk}$ for
$C\log(n)\le k\le c'\sqrt n$ follows by taking $\varepsilon=0$
in the bound obtained by Song and Wang
in~\cite[Theorem~1.5]{SongWangGap}.
The entries need not have the same distribution; the constants
and the lower bound on $n$ depend only on $T$.

For an $n\times n$ matrix $M$ with independent and identically
distributed non-constant subgaussian entries, the bound
$\Prob[\rank(M)\le n-k]\le e^{-c'nk}$ for $1\le k\le c\sqrt n$
was obtained by Rudelson in~\cite[Theorem~1.1]{Rudelson},
where $c,c'>0$ depend only on the distribution of the entries.
The restriction on $k$ was removed for Rademacher entries by
Hunter, Kwan, Sauermann, and Sawhney
in~\cite[Theorem~1.1]{HKSS}.
The same bound for all $1\le k\le n$, under a uniform bound below
one on the largest atom of each independent entry, was proved by
Polavarapu in~\cite[Theorem~1.1]{Polavarapu}.
These results require all entries to be independent.

For independent Bernoulli$(p)$ entries with fixed $p\in(0,1/2]$,
the sharp exponential rate for fixed $k$ was determined by
Jain, Sah, and Sawhney in~\cite[Theorem~1.1]{JSSRank}.
Song and Wang extended this result to $k=o(\sqrt{\log(n)})$
in~\cite[Theorem~1.1]{SongWang}.

\subsection{Proof ideas}\label{sec:proofideas}

Our main step compares column space probabilities using a common
auxiliary random subspace. This reduces the symmetric problem to
a corank bound for a matrix with independent entries.

For corank at least three, the principal submatrix argument of Han
in~\cite[Section~6]{Han} gives a large symmetric block with the same
rank as the original matrix and comparable corank.
The remaining columns, restricted to its rows, must lie in its column
space; they are independent of the block and of one another.
Lemma~\ref{lem:comparison} bounds this probability using the mass of
the column space under a fixed auxiliary distribution.
Small mass gives the required bound directly. We call the remaining
subspaces \emph{thick}.

Following the basis argument of Tao and Vu
in~\cite[Section~4]{TaoVu}, we sample independent vectors from this
distribution conditioned on lying in a fixed thick subspace,
and complete them to a basis with a small number of matrix columns.
We extend their indices to obtain an invertible principal block.
Fixing the corresponding columns leaves a large random symmetric
block and reduces the problem to a weighted sum over subspaces.
The weights are the probabilities that the sampled vectors and the
fixed columns span the given subspace.

Removing the conditioning from all sampled vectors can cost a factor
exponential in the square of the matrix size.
Instead, we construct an auxiliary matrix with the same fixed
columns and two off-diagonal blocks given by a rectangular matrix
with independent entries and its transpose, all other entries
are zero (see e.g., \cite{Han, CJMSLeast}). Comparing the entire kernel system bounds each conditional
column space probability by the probability that the auxiliary column
space is contained in the given subspace.
With the original columns fixed, the distribution of the auxiliary
matrix is the same for every subspace.

When the auxiliary column space is contained in the given subspace,
we remove conditioning only from the sampled vectors needed to span
that subspace together with the auxiliary column space.
For each choice of indices, the resulting spanning events are
disjoint for distinct subspaces. Lemma~\ref{lem:common-subspace}
applies the resulting summation bound to a common random subspace.
The number of vectors needed is controlled by the corank of the
auxiliary block with independent entries, which in turn must be
at least a constant multiple of the codimension of the given subspace.
Combining the rectangular corank bound
in~\cite[Corollary~1.2]{Polavarapu} with the bounds for the conditioning
loss and the number of choices gives the required probability bound.

For corank one, the relevant principal block $B$ is invertible, so
its column space gives no useful constraint on $Z$.
We use the Schur complement and the decoupling
inequality of Costello, Tao, and Vu \cite{CTV}
, in the reduction
developed by Ferber and Jain
\cite[Section~2.3]{FerberJain}
and by
Campos, Mattos, Morris, and Morrison
\cite[Lemma~2.1 and Appendix~A]{CMMM}
(see also \cite[Lemma~9.1 and Appendix~C]{CJMS})
.
Fixing one of the vectors and adding the corresponding fixed row
expresses the equation as the event that an independent random vector
lies in the augmented column space.
Apart from an exponentially small exceptional contribution, a normal
vector to this space has a nonzero coordinate among the row and column
indices of one of the chosen square blocks.
If the auxiliary column space is contained in this space, the
corresponding auxiliary block must be singular, allowing the same
comparison in Lemma~\ref{lem:augmented}. Corank two reduces to corank
one in a principal submatrix.

\subsection{Notation}

For a positive integer $n$, write $[n]=\{1,\ldots,n\}$,
and let $[0]=\varnothing$. For an $m\times n$ matrix $M$, write
$M_j$ for its $j$-th column. For $I\subseteq[m]$ and
$J\subseteq[n]$, let $M_{I\times J}$ be the submatrix with
row set $I$ and column set $J$. We write
$V_M:=\operatorname{span}(M_1,\ldots,M_n)$
for the column space of $M$, and let
$\corank(M)=\min\{m,n\}-\rank(M)$.
All vector spaces are over $\R$, and all vectors are column vectors.
For a subspace $V$, write
$V^\perp$ for its orthogonal complement.
For vectors, $\langle\cdot,\cdot\rangle$ is the usual inner
product and $\|\cdot\|_2$ is the corresponding norm.
The span of no vectors is $\{0\}$. We take the $0\times0$ matrix
to be invertible.

For a vector $x\in\R^m$ and $I\subseteq[m]$, write
$x_I=(x_i)_{i\in I}$, and let $\pi_I(x)=x_I$ be the coordinate
projection.

We use $\log$ for the natural logarithm. Constants and lower bounds
on dimensions associated with a lemma or proposition are indexed by
its number. These constants are fixed. Except where stated otherwise,
they depend only on $\beta$. We use $n_{\ref*{thm:main}}$ for the lower
bound on $n$ in the proof of Theorem~\ref{thm:main}. The parameter
$\varepsilon$ is chosen separately in Proposition~\ref{prop:thick}
and Lemma~\ref{lem:augmented}.

For a real random variable $X$, its largest atom is
$\sup_{x\in\R}\Prob[X=x]$.
For an event $E$, let $\mathbf1_E$ be its indicator.
We use $\Prob_U$ and $\E_U$ for probability and expectation
with respect to $U$.
For a probability distribution $\sigma$ on $\R^m$ and a subspace
$V\subseteq\R^m$ with $\sigma(V)>0$, write $\sigma(\,\cdot\mid V)$
for the law of a vector sampled from $\sigma$ conditioned on lying in $V$.
Thus, for every Borel set $E\subseteq\R^m$,
\[
 \sigma(E\mid V)=\frac{\sigma(E\cap V)}{\sigma(V)}.
\]

\subsection{Organization}

Section~\ref{sec:prelim} collects the preliminary tools.
Section~\ref{sec:reduction} reduces the case $k\ge3$ to a thick
column space bound, proved in Section~\ref{sec:thick}.
Lemma~\ref{lem:common-subspace} gives the general comparison;
Lemma~\ref{lem:conditional-bound} applies it to symmetric matrices
for any fixed probability distribution $\sigma$.
Section~\ref{sec:singularity} treats corank one using
Lemma~\ref{lem:augmented}, which allows any fixed nonzero added row,
and completes the proof.

\section{Preliminaries}\label{sec:prelim}

\subsection{Basic tools}

We first recall the bound of Odlyzko in~\cite{Odlyzko}, in the
form stated by Polavarapu in~\cite[Lemma~2.3]{Polavarapu}.

\begin{lemma}[Odlyzko bound]\label{lem:atom}
Let $N$ be a positive integer, let $0<\lambda<1$, and let
$\eta\in\R^N$ have independent coordinates, each with largest
atom at most $\lambda$. If $U\subseteq\R^N$ is a subspace of
codimension $\ell=N-\dim(U)$, then, for every $x\in\R^N$,
\[
 \Prob[\eta\in x+U]\le\lambda^\ell.
\]
\end{lemma}

The following rectangular rank bound was obtained by Polavarapu
in~\cite[Corollary~1.2]{Polavarapu}.
Its proof uses the determinant version of the theorem of
Hunter, Kwan, and Sauermann in~\cite[Theorem~1.2]{HKS}
for bounded corank.

\begin{lemma}[Rectangular corank bound]\label{lem:auxiliary-rank}
For every $0<\lambda<1$, there is a constant
$c_{\ref*{lem:auxiliary-rank}}=c_{\ref*{lem:auxiliary-rank}}(\lambda)>0$ such that the following holds.
Let $m\le n$ be positive integers,
and let $Z$ be an $m\times n$ matrix with independent real entries,
each with largest atom at most $\lambda$. For every integer
$1\le h\le m$,
\[
 \Prob[\rank(Z)\le m-h]\le e^{-c_{\ref*{lem:auxiliary-rank}}mh}.
\]
\end{lemma}

We use the following results on principal submatrices from the work of
Kodiyalam, Lam, and Swan
in~\cite[Theorem~2.9 and Corollary~2.10]{KLS}.

\begin{lemma}[Principal submatrices]\label{lem:principal}
Let $M$ be a real symmetric $n\times n$ matrix of rank $r$.
Then $M$ has an invertible principal submatrix of size $r$.
Moreover, if $I\subseteq[n]$ is the index set of $d\le r$ linearly
independent columns of $M$, there is a set $K$ with
$I\subseteq K\subseteq[n]$ and $|K|\le2d$
such that $M_{K\times K}$ is invertible.
\end{lemma}

We also use the following decoupling inequality of Costello, Tao,
and Vu in~\cite[Lemma~4.7]{CTV}.

\begin{lemma}[Decoupling]\label{lem:decoupling}
Let $u\in\R^r$ and $w\in\R^s$ be independent random vectors,
and let $u',w'$ be independent copies of $u,w$, respectively,
independent of $u,w$.
For every Borel set $E\subseteq\R^r\times\R^s$,
\[
 \Prob[(u,w)\in E]^4
 \le\Prob[(u,w),(u',w),(u,w'),(u',w')\in E].
\]
\end{lemma}

\subsection{A comparison for linear equations}

For a real random variable $\xi$, let $\xi'$ be an independent copy,
and let $\mu_\xi$ be the law of the product of $\xi-\xi'$ and an
independent Bernoulli$(1/2)$ random variable.

The following comparison uses the characteristic function argument
of Campos, Jenssen, Michelen, and Sahasrabudhe
in~\cite[Lemma~V.1 and proof of Lemma~III.7]{CJMSLeast}.
For earlier comparisons with lazy auxiliary variables (see e.g., \cite[Section~3]{TaoVu}).

\begin{lemma}[Comparison for linear equations]\label{lem:comparison}
Let $m,n$ be positive integers, and fix $L\in\R^{m\times n}$
and $y\in\R^m$. Let $X\in\R^n$ have independent real
coordinates, and let $\eta\in\R^n$ have independent coordinates,
where $\eta_j$ is distributed according to $\mu_{X_j}$ for each
$j\in[n]$. Then
\[
 \Prob[LX=y]\le\Prob[L\eta=0].
\]
\end{lemma}

\begin{proof}
\allowdisplaybreaks[1]
For each $j\in[n]$, let $\phi_j(t)=\E e^{itX_j}$, and let
$\psi_j$ be the characteristic function of $\eta_j$.
For every $t\in\R$,
\[
 \psi_j(t)=\frac{1+|\phi_j(t)|^2}{2}\ge|\phi_j(t)|.
\]
In particular, $\psi_j(t)$ is nonnegative.
Let $G$ be a standard Gaussian vector in $\R^m$, independent
of $X$ and $\eta$. For every $T>0$, the characteristic function
of $G$ and independence give
\begin{align*}
 \E\exp\left(-\frac{T^2}{2}\|LX-y\|_2^2\right)
 &=\E_G\left[e^{-iT\langle G,y\rangle}
       \prod_{j=1}^n\phi_j\bigl(T\langle L_j,G\rangle\bigr)\right]\\
 &\le\E_G\prod_{j=1}^n
       \left|\phi_j\bigl(T\langle L_j,G\rangle\bigr)\right|\\
 &\le\E_G\prod_{j=1}^n
       \psi_j\bigl(T\langle L_j,G\rangle\bigr)\\
 &=\E\exp\left(-\frac{T^2}{2}\|L\eta\|_2^2\right).
\end{align*}

The first inequality follows by taking absolute values.
Letting $T\to\infty$ gives the required probability bound,
since both random variables inside the expectations take values
in $[0,1]$.
\end{proof}

If $\xi$ has largest atom at most $\beta\in(0,1)$, then
$\Prob[\xi=\xi']\le\beta$ and $\mu_\xi(\{x\})\le1/2$ for every $x\ne0$.
Thus, with $\lambda=(1+\beta)/2$, we have
\begin{equation}\label{eq:atom-bound}
 \sup_{x\in\R}\mu_\xi(\{x\})
 =\mu_\xi(\{0\})=\frac{1+\Prob[\xi=\xi']}{2}\le\lambda<1.
\end{equation}
For Rademacher $\xi$, the probabilities of $0,2,-2$ under $\mu_\xi$
are $3/4,1/8,1/8$, respectively.

\section{Reduction to thick subspaces}\label{sec:reduction}

Throughout Sections~\ref{sec:reduction}--\ref{sec:singularity},
we work with fixed diagonal entries, and all bounds are uniform
in their values. We average over their joint distribution only at
the end of the proof of Theorem~\ref{thm:main}.

We first consider corank at least three, following the principal
submatrix argument of Han in~\cite[Section~6]{Han}.
Fix $\beta\in(0,1)$, let $\lambda=(1+\beta)/2$, and let $A$ satisfy
the hypotheses of Theorem~\ref{thm:main}.
For $3\le k\le n$, let
\[
 N=n-\lfloor k/2\rfloor,\qquad \ell=\lceil k/2\rceil.
\]
On $\{\corank(A)=k\}$, Lemma~\ref{lem:principal} gives an
invertible principal submatrix of size $n-k$.
Extend its index set to one of size $N$.
The corresponding principal submatrix $B$ satisfies
$n-k\le\rank(B)\le\rank(A)=n-k$, so $\corank(B)=\ell$.
This also holds when $n-k=0$, using the empty submatrix.

For each of the $\binom nN$ possible index sets, apply the same
permutation to the rows and columns and write
\[
 A=\begin{pmatrix}B&X\\X^T&\ast\end{pmatrix},
 \qquad X\in\R^{N\times(n-N)}.
\]
The blocks $B,X$ are independent.
The entries of $X$ are independent, each with largest atom at most
$\beta$.
If $\rank(A)=\rank(B)=n-k$, every column of $X$ lies in $V_B$.
For $j\in[n-N]$, let $\sigma_j$ be the product distribution obtained
from the distribution of $X_j$ by replacing each coordinate $\xi$
with an independent coordinate distributed according to $\mu_\xi$.
As in~\cite[Claim~4.4]{Polavarapu}, let
\[
 \sigma=\frac1{n-N}\sum_{j=1}^{n-N}\sigma_j.
\]
These distributions depend only on the entry distributions and
the fixed index set, and not on the realization of $B$.
By~\eqref{eq:atom-bound} and Lemma~\ref{lem:atom}, their average satisfies
\begin{equation}\label{eq:general-subspace}
 \sigma(U)\le\lambda^{N-\dim(U)}
 \quad\text{for every subspace }U\subseteq\R^N.
\end{equation}

For each fixed $B$, Lemma~\ref{lem:comparison} and the inequality
between arithmetic and geometric means give
\[
 \Prob[X_j\in V_B\text{ for every }j\in[n-N]\mid B]
 \le\prod_{j=1}^{n-N}\sigma_j(V_B)
 \le\sigma(V_B)^{n-N}.
\]
For $\varepsilon>0$, call a subspace $V\subseteq\R^N$ \emph{thin}
if $\sigma(V)<e^{-\varepsilon N}$, and \emph{thick} otherwise,
as in~\cite[Claims~2.2--2.3]{HKSS}.
These terms refer to the chosen distribution $\sigma$ and threshold.
On the thin event, the preceding conditional probability is at most
$e^{-\varepsilon N(n-N)}$. Splitting into the two cases gives
\begin{equation}\label{eq:thin-thick-reduction}
\begin{aligned}
 &\Prob[\corank(A)=k,\ \rank(B)=N-\ell]\\
 &\qquad\le e^{-\varepsilon N(n-N)}
   +\Prob[\rank(B)=N-\ell,\ \sigma(V_B)\ge e^{-\varepsilon N}].
\end{aligned}
\end{equation}
Since $N\ge n/2$ and $n-N\ge k/3$, the first term is at most
$e^{-\varepsilon nk/6}$. The remaining task is to bound the
probability that the symmetric block has a thick column space.

\begin{proposition}[Thick column space bound]\label{prop:thick}
Fix $\beta\in(0,1)$, and let $\lambda=(1+\beta)/2$.
There are constants $\varepsilon,c_{\ref*{prop:thick}}>0$ and a positive
integer $N_{\ref*{prop:thick}}$, depending only on $\beta$, such that the
following holds for every $N\ge N_{\ref*{prop:thick}}$.
Let $\sigma$ be a fixed probability distribution on $\R^N$
satisfying $\sigma(U)\le\lambda^{N-\dim(U)}$
for every subspace $U\subseteq\R^N$.
Let $B$ be an $N\times N$ matrix satisfying the hypotheses of
Theorem~\ref{thm:main} with fixed diagonal entries.
Then, for every integer $2\le\ell\le N$,
\begin{equation}\label{eq:general-thick}
 \Prob[\rank(B)=N-\ell,\ \sigma(V_B)\ge e^{-\varepsilon N}]
 \le e^{-c_{\ref*{prop:thick}}N\ell}.
\end{equation}
\end{proposition}

We prove Proposition~\ref{prop:thick} in Section~\ref{sec:thick}.
We first deduce the bound for corank at least three.

\begin{proposition}[Corank at least three]\label{prop:large-corank}
Fix $\beta\in(0,1)$. There are a constant
$c_{\ref*{prop:large-corank}}>0$ and a positive integer $n_{\ref*{prop:large-corank}}$,
depending only on $\beta$, such that the following holds.
Let $A$ satisfy the hypotheses of Theorem~\ref{thm:main} with fixed
diagonal entries. Then
\[
 \Prob[\rank(A)\le n-k]\le e^{-c_{\ref*{prop:large-corank}}nk}
\]
for every $n\ge n_{\ref*{prop:large-corank}}$ and every integer $3\le k\le n$.
\end{proposition}

\begin{proof}
Let $\varepsilon,c_{\ref*{prop:thick}},N_{\ref*{prop:thick}}$ be as in Proposition~\ref{prop:thick}.
Fix $3\le k\le n$, and take $N=n-\lfloor k/2\rfloor$ and
$\ell=\lceil k/2\rceil$.
Take $n$ sufficiently large that $N\ge n/2\ge N_{\ref*{prop:thick}}$.
For each principal index set, apply Proposition~\ref{prop:thick}
to the second term in~\eqref{eq:thin-thick-reduction}.
A union bound over the $\binom nN$ index sets gives
\begin{equation}\label{eq:exact-corank}
 \Prob[\corank(A)=k]
 \le\binom nN\left(
       e^{-\varepsilon N(n-N)}+e^{-c_{\ref*{prop:thick}}N\ell}\right).
\end{equation}
Take
\[
 c_{\ref*{prop:large-corank}}
 =\min\{\varepsilon/24,c_{\ref*{prop:thick}}/16\}>0.
\]
Since $N\ge n/2$, $n-N\ge k/3$, $\ell\ge k/2$, and
$\binom nN\le n^{n-N}\le n^{k/2}$, we obtain
\[
 \Prob[\corank(A)=k]
 \le2\exp\left(\frac{k}{2}\log(n)-4c_{\ref*{prop:large-corank}}nk\right)
 \le e^{-2c_{\ref*{prop:large-corank}}nk}
\]
when $\log(2n)\le 4c_{\ref*{prop:large-corank}}n$.
This lower bound on $n$ depends only on $\beta$ and is uniform in
$3\le k\le n$.
For sufficiently large $n$, we also have
\begin{align*}
 \Prob[\corank(A)\ge k]
 &\le\sum_{j=k}^{n}e^{-2c_{\ref*{prop:large-corank}}nj}\\
 &\le\frac{e^{-2c_{\ref*{prop:large-corank}}nk}}{1-e^{-2c_{\ref*{prop:large-corank}}n}}
 \le e^{-c_{\ref*{prop:large-corank}}nk},\qquad 3\le k\le n.
\end{align*}
Taking $n_{\ref*{prop:large-corank}}\ge2N_{\ref*{prop:thick}}$
sufficiently large proves the proposition.
\end{proof}

\section{Column space bounds}\label{sec:thick}

We prove Proposition~\ref{prop:thick} by reducing the thick event
to a weighted sum over column spaces and bounding it using a
common auxiliary subspace.

\subsection{Reduction to a weighted sum}\label{sec:spaces}

Let $\sigma$ be a fixed probability distribution on $\R^N$.
For each fixed subspace $V\subseteq\R^N$ with $\sigma(V)>0$ and
integer $t\ge0$, let $Y_1,\ldots,Y_t$ be independent vectors with
law $\sigma(\,\cdot\mid V)$, and write
\begin{equation}\label{eq:space-definition}
 W=\operatorname{span}(Y_1,\ldots,Y_t).
\end{equation}
For fixed $V$, the auxiliary vectors are independent of the random
matrix under consideration. When $V=V_B$ for a random matrix $B$,
we first fix $B$ and use this law for its column space.

The following lemma reduces the thick event to a weighted sum.
Its proof uses the basis argument of Tao and Vu
in~\cite[Section~4]{TaoVu}, as in~\cite[Claim~4.4]{Polavarapu},
and Lemma~\ref{lem:principal}.

\begin{lemma}[Reduction to a weighted sum]\label{lem:space-reduction}
Fix $0<\lambda<1$ and $0<\varepsilon\le-\log(\lambda)/128$.
There is a positive integer $N_{\ref*{lem:space-reduction}}$, depending
only on $\lambda,\varepsilon$, such that the following holds for
every $N\ge N_{\ref*{lem:space-reduction}}$.
Let $B$ be a real symmetric random $N\times N$ matrix, and let
$\sigma$ be a fixed probability distribution on $\R^N$ satisfying
\[
 \sigma(U)\le\lambda^{N-\dim(U)}
 \quad\text{for every subspace }U\subseteq\R^N.
\]
Fix $2\le\ell\le N$, let $r=N-\ell$, and assume that
\[
 \mathcal V_{\mathrm{thick}}
 =\{V\subseteq\R^N:\dim(V)=r,\ \sigma(V)\ge e^{-\varepsilon N}\}
\]
is nonempty. Let
\begin{equation}\label{eq:pool-dimensions}
 d=\left\lceil\frac{2\varepsilon N}{-\log(\lambda)}\right\rceil,
 \qquad t=r-d,
\end{equation}
and let $E=\{V_B\in\mathcal V_{\mathrm{thick}}\}$.
Then $0\le t<r$ and $2d\le\min\{r,N/4\}$.
For each $K\subseteq[N]$ with $d\le|K|\le2d$, let
$F=B_{[N]\times K}$. With $W$ as in~\eqref{eq:space-definition}
for each $V\in\mathcal V_{\mathrm{thick}}$, we have
\begin{equation}\label{eq:space-reduction}
 \Prob[E]\le
 2\sum_{\substack{K\subseteq[N]\\d\le|K|\le2d}}
 \E_F\Bigl[
   \mathbf1_{\{\det(B_{K\times K})\ne0\}}
   \sum_{\substack{V\in\mathcal V_{\mathrm{thick}}\\V_F\subseteq V}}
   \Prob_Y[V_F+W=V]\Prob_B[V_B=V\mid F]\Bigr].
\end{equation}
\end{lemma}

\begin{proof}
Since $\mathcal V_{\mathrm{thick}}$ is nonempty, the hypothesis on
$\sigma$ gives $e^{-\varepsilon N}\le\lambda^{N-r}$.
Hence $N-r\le\varepsilon N/(-\log(\lambda))$.
The definitions of $d,t$ then give $0\le t<r$ and
$2d\le\min\{r,N/4\}$ for sufficiently large $N$.

Now fix a realization of $B$ on $E$ and take auxiliary vectors with law
$\sigma(\,\cdot\mid V_B)$.
The first $j$ vectors span a subspace of $V_B$ of dimension at most
$j$. Thus
\begin{align*}
 &\Prob_Y\bigl[Y_{j+1}\in\operatorname{span}(Y_1,\ldots,Y_j)
       \mid Y_1,\ldots,Y_j\bigr]\\
 &\qquad=
 \frac{\sigma(\operatorname{span}(Y_1,\ldots,Y_j))}{\sigma(V_B)}
 \le e^{\varepsilon N}\lambda^{N-j}\\
 &\qquad\le e^{\varepsilon N}\lambda^d\le e^{-\varepsilon N},
 \qquad 0\le j<t.
\end{align*}
Here we used $\sigma(V_B)\ge e^{-\varepsilon N}$,
$N-j\ge d$ for $j<t=r-d\le N-d$, and
$\lambda^d\le e^{-2\varepsilon N}$ by~\eqref{eq:pool-dimensions}.
A union bound gives
\begin{equation}\label{eq:pool-success}
 \Prob_Y[\dim(W)=t]
 \ge1-te^{-\varepsilon N}\ge\frac12.
\end{equation}

On this event, choose $d$ columns of $B$ which, together with
$Y_1,\ldots,Y_t$, form a basis of $V_B$.
These columns are linearly independent. By Lemma~\ref{lem:principal},
their indices can be extended to a set $K$ with
$d\le|K|\le2d$ and $B_{K\times K}$ invertible.
Writing $F=B_{[N]\times K}$, we still have $V_F+W=V_B$
with all $t$ auxiliary vectors, since the added columns lie in $V_B$.
By~\eqref{eq:pool-success} and a union bound,
\begin{equation}\label{eq:basis-witness}
\begin{aligned}
 \frac12
 &\le\mathbb{P}_Y[\dim(W)=t]
 \le\mathbb{P}_Y[\exists K:\ V_F+W=V_B]\\
 &\le\sum_{\substack{K\subseteq[N]\\d\le|K|\le2d\\
                  \det(B_{K\times K})\ne0}}
 \mathbb{P}_Y[V_F+W=V_B].
\end{aligned}
\end{equation}
In the middle event, $K$ ranges over the index sets in the sum.

To average over $B$, we use the following observation, valid for any
fixed distribution $\sigma$ on $\R^N$ and integer $t\ge0$.
For a fixed subspace $S\subseteq\R^N$ and independent vectors
$\eta_1,\ldots,\eta_t$ with law $\sigma$, we have
\begin{equation}\label{eq:space-atom}
 \sigma(V)^t\Prob_Y[S+W=V]
 =\Prob_\eta\bigl[S+\operatorname{span}(\eta_1,\ldots,\eta_t)=V\bigr].
\end{equation}
Indeed, the event on the right requires all $t$ vectors to lie in
$V$.
For fixed $S,t$, these events are disjoint for distinct $V$,
so only countably many spaces with $\sigma(V)>0$ have
$\Prob_Y[S+W=V]>0$.
All sums over subspaces in this paper omit terms with zero
spanning probability.

Using the measurable structure on subspaces
in~\cite[Lemmas~2.1--2.2]{Polavarapu}, the probability on the
right-hand side of~\eqref{eq:space-atom} is measurable in $S,V$.
Dividing by $\sigma(V)^t>0$ shows that $\Prob_Y[S+W=V]$ is
measurable as well.
For fixed $F$, the sum is the integral of this measurable function
with respect to the conditional distribution of $V_B$.
Thus no measurable enumeration of the subspaces is needed.

The product $\mathbf1_E\Prob_Y[V_F+W=V_B]$ is defined as zero
outside $E$.
For each fixed $K$, let
$F=B_{[N]\times K}$. Multiply~\eqref{eq:basis-witness} by
$2\mathbf1_E$, take expectation over $B$, and condition on $F$.
The invertibility of $B_{K\times K}$ is determined by $F$.
Taking $N_{\ref*{lem:space-reduction}}$ sufficiently large
proves~\eqref{eq:space-reduction}.
\end{proof}

\begin{remark*}
For a matrix $B$ with all $N^2$ entries independent, fix
$K\subseteq[N]$, condition on $F=B_{[N]\times K}$, and let $V=V_F$.
Then
\[
 \Prob_B[V_B=V\mid F]
 =\prod_{j\notin K}\Prob_B[B_j\in V\mid F].
\]
For symmetric $B$, fixing columns also fixes rows, so the remaining
columns need not be independent. We use the spanning probabilities
as weights in~\eqref{eq:space-reduction} and bound the whole sum.
\end{remark*}

\subsection{A common subspace comparison}

To bound the inner sum in~\eqref{eq:space-reduction}, we remove
conditioning only from the auxiliary vectors needed to complete
a fixed subspace.
The following observation uses the counting argument of Tao and Vu
in~\cite[Section~4]{TaoVu} and Hunter, Kwan, Sauermann, and Sawhney
in~\cite[Claim~2.3]{HKSS}.

Let $\sigma$ be a probability distribution on $\R^m$ and
$S\subseteq\R^m$ a fixed subspace.
Fix integers $\dim(S)\le r\le m$ and $t\ge r-\dim(S)$.
For each $r$-dimensional subspace $V\supseteq S$ with $\sigma(V)>0$,
let $Y_1,\ldots,Y_t$ be independent vectors distributed according to
$\sigma(\,\cdot\mid V)$,
and let $W$ be their span, as in~\eqref{eq:space-definition}.

First fix $T\supseteq S$ with $\dim(T)\le r$, and let
$u=r-\dim(T)$, so $0\le u\le t$.
For each $V\supseteq T$, the event $S+W=V$ implies that some $u$
of the sampled vectors together with $T$ span $V$.
Let $\eta_1,\ldots,\eta_u$ be independent vectors with law $\sigma$.
A union bound over the $\binom tu$ choices gives
\[
 \sigma(V)^u\Prob_Y[S+W=V]
 \le\binom tu\Prob_\eta\bigl[
 T+\operatorname{span}(\eta_1,\ldots,\eta_u)=V\bigr].
\]
The factor $\sigma(V)^u$ removes the conditioning on the selected
vectors lying in $V$. The events on the right are disjoint for
distinct $V$, so
\begin{equation}\label{eq:space-sum}
 \sum_{\substack{V\subseteq\R^m,\ \dim(V)=r\\
                 T\subseteq V,\ \sigma(V)>0}}
 \sigma(V)^u\Prob_Y[S+W=V]\le\binom tu.
\end{equation}

We now apply this bound to a random subspace.
The same random subspace and event are used for every subspace
in the family.

\begin{lemma}[Common subspace comparison]\label{lem:common-subspace}
In the above setting, fix a positive integer $N$ and $\varepsilon>0$.
Let $\mathcal V$ be a family of $r$-dimensional subspaces
$V\supseteq S$ with $\sigma(V)\ge e^{-\varepsilon N}$.
Let $T\supseteq S$ be a random subspace of $\R^m$, and let
$\mathcal E$ be an event on the same probability space.
Writing $u=r-\dim(T)$, we have
\begin{equation}\label{eq:common-subspace}
 \sum_{V\in\mathcal V}\Prob_Y[S+W=V]\Prob[\mathcal E,\ T\subseteq V]
 \le\E\left[
 \mathbf1_{\mathcal E\cap\{u\ge0\}}\binom tu e^{\varepsilon Nu}\right],
\end{equation}
where the integrand is defined as zero when $u<0$.
\end{lemma}

\begin{proof}
The left-hand side of~\eqref{eq:common-subspace} equals
\[
 \E\left[\mathbf1_{\mathcal E}
 \sum_{\substack{V\in\mathcal V\\T\subseteq V}}
 \Prob_Y[S+W=V]\right].
\]
For each realization of $T$, the inner sum is zero if $\dim(T)>r$.
Otherwise, $0\le u\le r-\dim(S)\le t$, and
$\sigma(V)^{-u}\le e^{\varepsilon Nu}$.
Equation~\eqref{eq:space-sum} therefore bounds the sum by
$\binom tu e^{\varepsilon Nu}$, proving~\eqref{eq:common-subspace}.
\end{proof}

\subsection{A conditional column space bound}\label{sec:conditional}

We now bound the inner sum in~\eqref{eq:space-reduction}.

For the following lemma, fix $\beta\in(0,1)$.
Let $B$ be an $N\times N$ matrix satisfying the hypotheses of
Theorem~\ref{thm:main} with fixed diagonal entries, and let $\sigma$
be any fixed probability distribution on $\R^N$.
Fix $K\subseteq[N]$ with $|K|\le N/4$, and condition on
$F=B_{[N]\times K}$ with $B_{K\times K}$ invertible.
Under this conditional law, the entries determined by $F$ are fixed,
and the remaining entries above the diagonal have their original
independent distributions. We use this law for every value of $F$
consistent with symmetry and the fixed diagonal entries, including
values of probability zero.
For $\varepsilon>0$ and integers $2\le\ell\le N$ and $0\le t\le N$,
let $\mathcal V$ be the family of $(N-\ell)$-dimensional subspaces
$V\supseteq V_F$ with $\sigma(V)\ge e^{-\varepsilon N}$.
For each $V\in\mathcal V$, let $W$ be as in~\eqref{eq:space-definition}.

\begin{lemma}[Weighted conditional column space bound]\label{lem:conditional-bound}
Under the above assumptions, there are a constant
$c_{\ref*{lem:conditional-bound}}>0$ and a positive integer
$N_{\ref*{lem:conditional-bound}}$, depending only on $\beta$, such that,
for every $N\ge N_{\ref*{lem:conditional-bound}}$ and
every $0<\varepsilon\le c_{\ref*{lem:conditional-bound}}$,
\begin{equation}\label{eq:conditional-bound}
 \sum_{V\in\mathcal V}
 \Prob_Y[V_F+W=V]\Prob_B[V_B=V\mid F]\le e^{-c_{\ref*{lem:conditional-bound}}N\ell}.
\end{equation}
\end{lemma}

\begin{proof}
Let $\lambda=(1+\beta)/2$ and $c_{\ref*{lem:auxiliary-rank}}=c_{\ref*{lem:auxiliary-rank}}(\lambda)$ be as in
Lemma~\ref{lem:auxiliary-rank}. Take $c_{\ref*{lem:conditional-bound}}=c_{\ref*{lem:auxiliary-rank}}/48$.
Fix $K,F$. Since $\dim(V_F)=|K|$ and $\dim(W)\le t$, we may
assume $|K|\le N-\ell\le|K|+t$, as otherwise the sum is zero.
All probabilities involving $B$ below use this conditional law.

Choose a fixed partition $[N]\setminus K=R_1\sqcup R_2$ with
$|R_1|\le|R_2|\le|R_1|+1$, and let $D=B_{R_1\times R_2}$.
Condition further on the remaining entries above the diagonal
outside $D$ and $D^T$.
The entries of $D$ remain independent with their original distributions.

Let $Z$ be a matrix of the same size as $D$, independent of $B$,
with independent entries $Z_{ij}$ distributed according to $\mu_{D_{ij}}$.
In the coordinate order $K,R_1,R_2$, write $B$ and define $H$ by
\[
 B=\begin{pmatrix}
  B_{K\times K}&*&*\\
  B_{R_1\times K}&*&D\\
  B_{R_2\times K}&D^T&*
 \end{pmatrix},
 \qquad
 H=\begin{pmatrix}
  B_{K\times K}&0&0\\
  B_{R_1\times K}&0&Z\\
  B_{R_2\times K}&Z^T&0
 \end{pmatrix}\in\R^{N\times N}.
\]
The first block column of both matrices is $F$, and the blocks
denoted by $*$ are fixed. The symmetric block formed by $Z$ and $Z^T$
is of the form
used by Campos, Jenssen, Michelen, and Sahasrabudhe
in~\cite[Sections~II.4 and~V]{CJMSLeast} and, for Rademacher
entries, in~\cite[Appendix~B]{CJMS}; see also~\cite[Section~4.2]{Han}.

Fix $V\in\mathcal V$, and let the columns of
$Q\in\R^{N\times\ell}$ form a basis of $V^\perp$.
Write $Q_1,Q_2$ for the row blocks of $Q$ indexed by $R_1,R_2$.
The event $V_B=V$ implies $B^TQ=0$.
Since $V_F\subseteq V$, the equations $F^TQ=0$ already hold.
Under the further conditioning, move the fixed terms in the remaining
equations to the right-hand side, obtaining a fixed matrix $M$.
This is a linear system in the independent entries of $D$, with
fixed coefficients determined by $Q_1,Q_2$.
Lemma~\ref{lem:comparison} therefore gives
\[
 \Prob_D[B^TQ=0]
 =\Prob_D\left[\begin{pmatrix}DQ_2\\D^TQ_1\end{pmatrix}=M\right]
 \le\Prob_Z\left[\begin{pmatrix}ZQ_2\\Z^TQ_1\end{pmatrix}=0\right].
\]
Since $F^TQ=0$, the event on the right is exactly $H^TQ=0$,
equivalently $V_H\subseteq V$.
With $F$ fixed, the law of $H$ is the same for every $V\in\mathcal V$
and does not depend on the further conditioning.
Averaging over the entries fixed by the further conditioning gives
\begin{equation}\label{eq:common-comparison}
 \Prob_B[V_B=V]\le\Prob_H[V_H\subseteq V].
\end{equation}
Since $B_{K\times K}$ is invertible and the columns of $H$ with
indices outside $K$ are zero on the rows indexed by $K$,
we have $\rank(H)=|K|+2\rank(Z)$.

For fixed $H$, let $u=N-\ell-\rank(H)$ and
$h=\corank(Z)=|R_1|-\rank(Z)$. Then
\[
 u=(|R_2|-|R_1|)+2h-\ell.
\]
If $u<0$, no space in $\mathcal V$ contains $V_H$.
Otherwise, $0\le u\le N-\ell-|K|\le t$.
Since $|R_2|-|R_1|\in\{0,1\}$ and $\ell\ge2$, on $\{u\ge0\}$ we have
\[
 u\le2h,\qquad
 h\ge\lfloor\ell/2\rfloor\ge\ell/3.
\]
Apply Lemma~\ref{lem:common-subspace} to~\eqref{eq:common-comparison}
with $m=N$, $S=V_F$, $T=V_H$, $r=N-\ell$, and
$\mathcal E$ the whole probability space. It gives
\begin{equation}\label{eq:common-space-sum}
\begin{aligned}
 &\sum_{V\in\mathcal V}
   \Prob_Y[V_F+W=V]\Prob_B[V_B=V]\\
 &\qquad\le\E_H\left[
   \mathbf1_{\{u\ge0\}}\binom tu e^{\varepsilon Nu}\right].
\end{aligned}
\end{equation}
The expression inside the expectation is defined as zero when
$u<0$.
For every nonzero term, both the number of ways to choose the
auxiliary vectors and the loss from removing conditioning
are bounded in terms of $h$:
\[
 \binom tu\le N^u\le N^{2h},\qquad
 e^{\varepsilon Nu}\le e^{2\varepsilon Nh}.
\]
Every entry of $Z$ has largest atom at most $\lambda$
by~\eqref{eq:atom-bound}.
Lemma~\ref{lem:auxiliary-rank} gives
\begin{align*}
 &\E_H\left[
   \mathbf1_{\{u\ge0\}}\binom tu e^{\varepsilon Nu}
 \right]\\
 &\qquad\le
 \sum_{h=\lfloor\ell/2\rfloor}^{|R_1|}
   N^{2h}e^{2\varepsilon Nh}\Prob[\corank(Z)=h]\\
 &\qquad\le
 \sum_{h=\lfloor\ell/2\rfloor}^{|R_1|}e^{-c_{\ref*{lem:auxiliary-rank}}|R_1|h/2}
 \le2e^{-c_{\ref*{lem:auxiliary-rank}}N\ell/18}
 \le e^{-c_{\ref*{lem:conditional-bound}}N\ell}
\end{align*}
for sufficiently large $N$.
Here we used $|R_1|\ge N/3$, $\lfloor\ell/2\rfloor\ge\ell/3$, and
$2\log(N)+2\varepsilon N\le c_{\ref*{lem:auxiliary-rank}}|R_1|/2$.
Taking $N_{\ref*{lem:conditional-bound}}$ sufficiently large
proves~\eqref{eq:conditional-bound}.
\end{proof}

\subsection{Proof of the thick column space bound}

\begin{proof}[Proof of Proposition~\ref{prop:thick}]
Let $\lambda=(1+\beta)/2$, and let $c_{\ref*{lem:conditional-bound}}$
be as in Lemma~\ref{lem:conditional-bound}.
Choose
\begin{equation}\label{eq:epsilon}
 0<\varepsilon\le
 \min\{-\log(\lambda)/128,c_{\ref*{lem:conditional-bound}}\}
\end{equation}
sufficiently small that, with $d$ as in~\eqref{eq:pool-dimensions},
\begin{equation}\label{eq:small-index-count}
 2\sum_{s=d}^{2d}\binom Ns
 \le2(d+1)\left(\frac{eN}{2d}\right)^{2d}
 \le e^{c_{\ref*{lem:conditional-bound}}N/2}
\end{equation}
for sufficiently large $N$.
Such a choice is possible since $d/N\to2\varepsilon/(-\log(\lambda))$
as $N\to\infty$ and $x\log(e/x)\to0$ as $x\to0$.
Choose
\[
 N_{\ref*{prop:thick}}\ge
 \max\{N_{\ref*{lem:space-reduction}},N_{\ref*{lem:conditional-bound}}\}
\]
sufficiently large that~\eqref{eq:small-index-count} holds for
every $N\ge N_{\ref*{prop:thick}}$.

Fix $N\ge N_{\ref*{prop:thick}}$ and $2\le\ell\le N$, and let $r=N-\ell$.
Let $\mathcal V_{\mathrm{thick}}$ be the family of $r$-dimensional
subspaces $V\subseteq\R^N$ with
$\sigma(V)\ge e^{-\varepsilon N}$.
If this family is empty, the result is immediate.
Otherwise, let $d,t$ be as in~\eqref{eq:pool-dimensions}.
Let $E=\{V_B\in\mathcal V_{\mathrm{thick}}\}$.
Apply Lemma~\ref{lem:space-reduction} with these parameters.
For each $K$ in~\eqref{eq:space-reduction}, we have
$|K|\le2d\le N/4$ and $0\le t\le N$.
Lemma~\ref{lem:conditional-bound} bounds the inner sum by
$e^{-c_{\ref*{lem:conditional-bound}}N\ell}$, uniformly in $F$ whenever
$B_{K\times K}$ is invertible.
Consequently,
\begin{align*}
 \Prob[E]
 &\le2\sum_{s=d}^{2d}\binom Ns e^{-c_{\ref*{lem:conditional-bound}}N\ell}\\
 &\le e^{c_{\ref*{lem:conditional-bound}}N/2-c_{\ref*{lem:conditional-bound}}N\ell}
 \le e^{-c_{\ref*{lem:conditional-bound}}N\ell/2}.
\end{align*}
Taking $c_{\ref*{prop:thick}}=c_{\ref*{lem:conditional-bound}}/2$ proves the proposition.
\end{proof}

\section{Corank one and completion of the proof}\label{sec:singularity}

For corank one, an invertible principal block leaves only one
column outside it. Its column space is the whole ambient space, so the
remaining column always lies in it. This condition therefore gives no
probability bound.
Adding the transpose of that column below the block gives a hyperplane
as the new column space. However, the added row and the vector
required to lie in this space depend on the same random column.
Fixing the added row would also fix that vector.

The Schur complement gives a quadratic equation in the remaining
column. We split the column into a small and a large coordinate
block and apply decoupling, reducing the problem to a bilinear
equation in two independent difference vectors;
related reductions are used in~\cite[Sections~5 and~8.2]{Vershynin}.
Fixing a nonzero value of the difference vector on the small block
determines the added row. The bilinear equation then says that
the random vector from the large block, extended by zeros, lies
in the column space of the augmented matrix. This vector is
independent of that matrix.

The law of this vector is used in both parts of the thin/thick
argument. If it assigns small mass to the column space, the
probability that the vector lies in that space is small.
Lemma~\ref{lem:augmented} bounds the probability that this mass is large.
We keep most coordinates in the large block so that the basis
argument needs only a small number of matrix columns.
The small block still contains a fixed positive fraction of the
coordinates, so the probability that its difference vector is zero
is exponentially small.

\subsection{Corank one}

\begin{lemma}[Corank one]\label{lem:corank-one}
Fix $\beta\in(0,1)$. There are a constant
$c_{\ref*{lem:corank-one}}>0$ and a positive integer $n_{\ref*{lem:corank-one}}$,
depending only on $\beta$, such that the following holds.
Let $A$ be an $n\times n$ matrix satisfying the hypotheses of
Theorem~\ref{thm:main} with fixed diagonal entries. Then
\[
 \Prob[\rank(A)=n-1]\le e^{-c_{\ref*{lem:corank-one}}n},\qquad n\ge n_{\ref*{lem:corank-one}}.
\]
\end{lemma}

\begin{proof}
If $\corank(A)=1$, Lemma~\ref{lem:principal} gives an invertible
principal submatrix of size $n-1$.
There are at most $n$ choices of its index set. Fix one and write
$A=\left(\begin{smallmatrix}B&x\\x^T&a\end{smallmatrix}\right)$,
with $N=n-1$. The vector $x$ has independent coordinates and is
independent of $B$, and $a$ is fixed.
On $\{\det(B)\ne0\}$, the Schur complement identity gives
$\det(A)=0$ if and only if $x^TB^{-1}x=a$.
Let $\varepsilon,c_{\ref*{lem:augmented}},N_{\ref*{lem:augmented}}$ be as in Lemma~\ref{lem:augmented} below.
For $N\ge N_{\ref*{lem:augmented}}$, split $[N]=I\sqcup J$ with $|J|=\lfloor\varepsilon N\rfloor$.
In the coordinate order $I,J$, write $x=(x_I,x_J)$.
Let $x'$ be a copy of $x$, independent of $B,x$, and let
$y=x-x'$.

For fixed invertible $B$, apply Lemma~\ref{lem:decoupling}
to the independent coordinate subvectors $x_I,x_J$ and the event
$x^TB^{-1}x=a$.
On the resulting intersection of four events, add the equations
for $x,x'$ and subtract those for $(x'_I,x_J)$ and $(x_I,x'_J)$.
The terms involving only one coordinate block cancel, giving
$2y_I^T(B^{-1})_{I\times J}y_J=0$. Hence
\[
 \Prob_x[x^TB^{-1}x=a]^4
 \le\Prob_y[y_I^T(B^{-1})_{I\times J}y_J=0].
\]
Averaging over $B,y_J$ and applying H\"older's inequality gives
\begin{equation}\label{eq:schur-decoupling}
 \Prob[\det(B)\ne0,\ x^TB^{-1}x=a]^4
 \le\E_{B,y_J}\left[
 \mathbf1_{\{\det(B)\ne0\}}
 \Prob_{y_I}[y_I^T(B^{-1})_{I\times J}y_J=0]\right].
\end{equation}
The expression inside the expectation is defined as zero when
$B$ is singular.
Conditioning on $x'$ shows that taking differences does not
increase the largest atom. Thus $y$ has independent coordinates,
each with largest atom at most $\beta$, and is independent of $B$.
In particular, $\Prob[y_J=0]\le\beta^{|J|}$.

Fix a nonzero value of $y_J$, and let $v=(0_I,y_J)$,
$C=\left(\begin{smallmatrix}B\\v^T\end{smallmatrix}\right)$,
and $\sigma$ be the law of $(y_I,0_J,0)$.
With $y_J$ fixed, neither $v$ nor $\sigma$ depends on $B$.
On $\{\det(B)\ne0\}$,
$V_C^\perp=\operatorname{span}((B^{-1}v,-1))$, so
\[
 \sigma(V_C)=\Prob_{y_I}[y_I^T(B^{-1})_{I\times J}y_J=0].
\]
We verify the hypothesis on $\sigma$ in Lemma~\ref{lem:augmented}.
Let $p=|I|=N-\lfloor\varepsilon N\rfloor$.
For every subspace $U\subseteq\R^{N+1}$, the set
$\{z\in\R^p:(z,0_J,0)\in U\}$
is a subspace of dimension at most $\dim(U)$.
Applying Lemma~\ref{lem:atom} to $y_I$ gives
\[
 \sigma(U)\le\beta^{p-\dim(U)}.
\]
On $\{\det(B)\ne0\}$, split according to whether
$\sigma(V_C)\ge e^{-\varepsilon N}$.
Lemma~\ref{lem:augmented} gives
\[
 \E_B[\mathbf1_{\{\det(B)\ne0\}}\sigma(V_C)]
 \le e^{-\varepsilon N}+e^{-c_{\ref*{lem:augmented}}N}
\]
uniformly in the fixed nonzero value of $y_J$.
Set
\[
 c_{\ref*{lem:corank-one}}
 =\frac1{16}\min\{-\varepsilon\log(\beta),
                   \varepsilon,c_{\ref*{lem:augmented}}\}>0.
\]
The right-hand side of~\eqref{eq:schur-decoupling} is at most
\[
 \beta^{\lfloor\varepsilon N\rfloor}
 +e^{-\varepsilon N}+e^{-c_{\ref*{lem:augmented}}N}
 \le e^{-8c_{\ref*{lem:corank-one}}N}
\]
for sufficiently large $N=n-1$.
Taking fourth roots and summing over the $n$ principal index sets gives
\[
 \Prob[\corank(A)=1]
 \le n e^{-2c_{\ref*{lem:corank-one}}(n-1)}
 \le e^{-c_{\ref*{lem:corank-one}}n}
\]
for sufficiently large $n$.
Taking $n_{\ref*{lem:corank-one}}\ge N_{\ref*{lem:augmented}}+1$
sufficiently large proves the lemma.
\end{proof}

\subsection{An augmented column space bound}\label{sec:augmented}

Fix $\beta\in(0,1)$, let $B$ be an $N\times N$ matrix satisfying
the hypotheses of Theorem~\ref{thm:main} with fixed diagonal entries,
and let $v\in\R^N\setminus\{0\}$ be fixed.
Write $C=\left(\begin{smallmatrix}B\\v^T\end{smallmatrix}\right)$.

\begin{lemma}[Augmented column space bound]\label{lem:augmented}
There are $0<\varepsilon<1/2$, $c_{\ref*{lem:augmented}}>0$, and a positive
integer $N_{\ref*{lem:augmented}}$, depending only on $\beta$, such that, for
$N\ge N_{\ref*{lem:augmented}}$ and every fixed
probability distribution $\sigma$ on $\R^{N+1}$ satisfying
\[
 \sigma(U)\le\beta^{p-\dim(U)},\qquad
 p=N-\lfloor\varepsilon N\rfloor,
\]
for every subspace $U\subseteq\R^{N+1}$, we have
\begin{equation}\label{eq:augmented-thick}
 \Prob[\det(B)\ne0,\ \sigma(V_C)\ge e^{-\varepsilon N}]
 \le e^{-c_{\ref*{lem:augmented}}N}.
\end{equation}
\end{lemma}

\begin{proof}
Let $c_{\ref*{lem:auxiliary-rank}}$ be the constant in
Lemma~\ref{lem:auxiliary-rank} for $\lambda=(1+\beta)/2$, and take
\[
 c_{\ref*{lem:augmented}}
 =\frac14\min\{-3\log(\beta)/4,c_{\ref*{lem:auxiliary-rank}}/8\}>0.
\]
For $\varepsilon>0$, let
\[
 d=\lfloor\varepsilon N\rfloor
   +\left\lceil\frac{2\varepsilon N}{-\log(\beta)}\right\rceil.
\]
By the binomial bound in the proof of Proposition~\ref{prop:thick},
choose $0<\varepsilon<1/2$ sufficiently small, with
$\varepsilon\le c_{\ref*{lem:auxiliary-rank}}/48$, so that
\[
 2d\le N/4,\qquad
 2\sum_{s=d}^{2d}\binom Ns\le e^{2c_{\ref*{lem:augmented}}N}
\]
for all sufficiently large $N$.
Let
\[
 \mathcal V
 =\{V\subseteq\R^{N+1}:\dim(V)=N,\
       \sigma(V)\ge e^{-\varepsilon N}\},
\]
and let $E=\{\det(B)\ne0,\ V_C\in\mathcal V\}$.
Let $t=N-d$. For each $V\in\mathcal V$, use the auxiliary vectors
and $W$ from Section~\ref{sec:spaces}, now in $\R^{N+1}$.

Fix $B$ on $E$ and sample the auxiliary vectors from
$\sigma(\,\cdot\mid V_C)$. Conditional on the first $j$ vectors,
the probability that the next lies in their span is at most
\[
 e^{\varepsilon N}\beta^{p-j}\le e^{-\varepsilon N},
 \qquad 0\le j<t.
\]
As in~\eqref{eq:pool-success}, the $t$ auxiliary vectors are
linearly independent with probability at least $1/2$.
On this event, complete them to a basis of $V_C$ with $d$ columns of $C$.
Since $B$ is invertible, Lemma~\ref{lem:principal}, applied to $B$,
extends their indices to $K$ with $d\le|K|\le2d$ and
$B_{K\times K}$ invertible.
For each fixed $K$, let $F=C_{[N+1]\times K}$.
The countability and measurability argument in the proof of
Lemma~\ref{lem:space-reduction} allows us to average over $B$
and condition on $F$, giving
\begin{equation}\label{eq:augmented-reduction}
 \Prob[E]\le
 2\sum_{\substack{K\subseteq[N]\\d\le|K|\le2d}}
 \E_F\Bigl[\mathbf1_{\{\det(B_{K\times K})\ne0\}}
 \sum_{\substack{V\in\mathcal V\\V_F\subseteq V}}
 \Prob_Y[V_F+W=V]\Prob_B[E,\ V_C=V\mid F]\Bigr].
\end{equation}

Fix $K$, let $s=|K|$, $R=[N]\setminus K$, and
$G=B_{[N]\times K}$.
Since $F=\left(\begin{smallmatrix}G\\v_K^T\end{smallmatrix}\right)$
and $v$ is fixed, conditioning on $F$ is equivalent to conditioning
on $G$. On $E$,
\[
 V_C^\perp=\operatorname{span}((B^{-1}v,-1)),\qquad
 (B^{-1}v)_R=0\quad\Longleftrightarrow\quad v\in V_G.
\]
We first bound the contribution to~\eqref{eq:augmented-reduction}
from values of $F$ for which $v\in V_G$.
Under the original law of $B$, condition only on a fixed invertible
value of $B_{K\times K}$. If $v\in V_G$, its coefficients in the columns
of $G$ must be $B_{K\times K}^{-1}v_K$.
If $v_K=0$, then $v\in V_G$ is impossible since $v\ne0$.
Otherwise the coefficient vector is fixed and nonzero.
Its inner products with the rows of $B_{R\times K}$ are independent,
each with largest atom at most $\beta$ by Lemma~\ref{lem:atom}.
Averaging over $B_{K\times K}$ gives
\[
 \Prob_F[\det(B_{K\times K})\ne0,\ v\in V_G]\le\beta^{N-s}.
\]
The inner sum in~\eqref{eq:augmented-reduction} is at most
$\Prob_B[E\mid F]\le1$.

Now fix $F$ with $B_{K\times K}$ invertible and $v\notin V_G$.
Use the conditional law from the proof of Lemma~\ref{lem:conditional-bound}.
The entries of $B$ above the diagonal not determined by $F$
remain independent with their original distributions.
It suffices to consider $V\in\mathcal V$ with $V_F\subseteq V$ and
$\Prob_B[E,\ V_C=V\mid F]>0$.
The preceding equivalence implies $\pi_R(V^\perp)\ne\{0\}$
for each such space.
Let $m_0=\lfloor |R|/2\rfloor$.
Choose at most two fixed subsets of $R$, each of size $2m_0$,
whose union is $R$. Assign each remaining space to the first
of these sets on which its normal space has nonzero projection.
Fix one resulting family, write $R_0$ for its associated set,
and fix a partition $R_0=P\sqcup Q$ with $|P|=|Q|=m_0$.
Let $D=B_{P\times Q}$.
As in the proof of Lemma~\ref{lem:conditional-bound}, condition
further on the remaining entries above the diagonal outside $D$ and $D^T$.
Let $Z\in\R^{m_0\times m_0}$ have independent entries with
$Z_{ij}$ distributed according to $\mu_{D_{ij}}$, independently of $B$.
Define $H\in\R^{(N+1)\times N}$ by setting
$H_{[N+1]\times K}=F$ and replacing the indicated block of $C$ as follows:
\[
 C_{R_0\times R_0}=
 \begin{pmatrix}*&D\\D^T&*\end{pmatrix},
 \qquad
 H_{R_0\times R_0}=
 \begin{pmatrix}0&Z\\Z^T&0\end{pmatrix}.
\]
The blocks of $C$ denoted by $*$ are fixed under this conditioning.
The displayed submatrices use the coordinate order $P,Q$.
All other entries of $H$ are zero; in particular, its last row
is $v_K^T$ on $K$ and zero elsewhere.
For fixed $V$ in this family, let $w$ span $V^\perp$.
Since $V_F\subseteq V$, we have $F^Tw=0$.
The event $E\cap\{V_C=V\}$ implies
$C^Tw=Bw_{[N]}+w_{N+1}v=0$.
Move all fixed terms, including $w_{N+1}v$, to the right-hand side
and apply Lemma~\ref{lem:comparison} once to the whole system
in the entries of $D$. This gives
\[
 \begin{aligned}
 \Prob_D[C^Tw=0]
 &\le\Prob_Z[Zw_{Q}=0,\ Z^Tw_{P}=0]\\
 &=\Prob_H[V_H\subseteq V],
 \end{aligned}
\]
where the last equality uses $F^Tw=0$.
With $F$ fixed, the law of $H$ is common to all spaces in this family
and does not depend on the further conditioning.
Averaging over the entries fixed by the further conditioning gives
\[
 \Prob_B[E,\ V_C=V\mid F]\le\Prob_H[V_H\subseteq V].
\]
If $V_H\subseteq V$ and $Z$ were invertible, then $H^Tw=0$ would
force $w_{R_0}=0$, contrary to $\pi_{R_0}(V^\perp)\ne\{0\}$.
Thus $V_H\subseteq V$ for any space in this family forces
$\corank(Z)\ge1$.

Since $B_{K\times K}$ is invertible and the columns of $H$ with
indices outside $K$ are zero on the rows indexed by $K$,
we have $\rank(H)=s+2\rank(Z)$.
Let $h=\corank(Z)$ and
$u=N-\rank(H)=|R|-2m_0+2h$.
Thus $0\le u\le N-s\le t$ and $u\le2h+1$.
Apply Lemma~\ref{lem:common-subspace} with $m=N+1$,
$S=V_F$, $T=V_H$, $r=N$, and
$\mathcal E=\{h\ge1\}$.
The preceding comparison and~\eqref{eq:atom-bound}, together with
Lemma~\ref{lem:auxiliary-rank}, bound the sum for this family by
\begin{align*}
 \E_H\left[\mathbf1_{\{h\ge1\}}\binom tu e^{\varepsilon Nu}\right]
 &\le Ne^{\varepsilon N}\sum_{h=1}^{m_0}
       e^{(2\log(N)+2\varepsilon N-c_{\ref*{lem:auxiliary-rank}}m_0)h}\\
 &\le e^{-c_{\ref*{lem:auxiliary-rank}}N/8},
\end{align*}
since $m_0\ge N/3$, $\varepsilon\le c_{\ref*{lem:auxiliary-rank}}/48$, and $N$ is sufficiently large.
Including the exceptional values of $F$ and the at most two
families in~\eqref{eq:augmented-reduction}, the index count above gives
\[
 \Prob[E]\le
 e^{2c_{\ref*{lem:augmented}}N}
 \left(\beta^{3N/4}+2e^{-c_{\ref*{lem:auxiliary-rank}}N/8}\right)
 \le3e^{-2c_{\ref*{lem:augmented}}N}
 \le e^{-c_{\ref*{lem:augmented}}N}.
\]
Taking $N_{\ref*{lem:augmented}}$ sufficiently large proves the lemma.
\end{proof}

\subsection{Completion of the proof}

\begin{proof}[Proof of Theorem~\ref{thm:main}]
First suppose that the diagonal entries are fixed.
Let $c_{\ref*{lem:corank-one}},n_{\ref*{lem:corank-one}}$ be as in Lemma~\ref{lem:corank-one}, and let
$c_{\ref*{prop:large-corank}},n_{\ref*{prop:large-corank}}$ be as in Proposition~\ref{prop:large-corank}.
Let $c=\min\{c_{\ref*{lem:corank-one}}/4,c_{\ref*{prop:large-corank}}\}>0$.
Choose $n_{\ref*{thm:main}}\ge\max\{n_{\ref*{lem:corank-one}}+1,n_{\ref*{prop:large-corank}},3\}$ sufficiently large, depending
only on $\beta$, and first suppose that $n\ge n_{\ref*{thm:main}}$.

On $\{\corank(A)=2\}$, Lemma~\ref{lem:principal} gives an
invertible principal submatrix of size $n-2$.
Adding one index gives a principal submatrix of size $n-1$
and rank $n-2$.
A union bound over the $n$ possible index sets and
Lemma~\ref{lem:corank-one} give
\[
 \Prob[\corank(A)=2]
 \le n e^{-c_{\ref*{lem:corank-one}}(n-1)}.
\]
For sufficiently large $n_{\ref*{thm:main}}$, we obtain
\begin{align*}
 \Prob[\det(A)=0]
 &\le e^{-c_{\ref*{lem:corank-one}}n}
       +n e^{-c_{\ref*{lem:corank-one}}(n-1)}+e^{-3c_{\ref*{prop:large-corank}}n}
 \le e^{-2cn}.
\end{align*}
Consequently, for $k=1,2$,
\[
 \Prob[\corank(A)\ge k]\le e^{-2cn}\le e^{-cnk}.
\]
For $3\le k\le n$, Proposition~\ref{prop:large-corank} gives
the same bound since $c\le c_{\ref*{prop:large-corank}}$.

Choose $0<c_\beta\le\min\{c,\log(2)/n_{\ref*{thm:main}}^2\}$.
The required bound holds for $n\ge n_{\ref*{thm:main}}$.
For $n<n_{\ref*{thm:main}}$ and $1\le k\le n$, we have
$c_\beta nk\le c_\beta n_{\ref*{thm:main}}^2\le\log(2)$, so
$2e^{-c_\beta nk}\ge1$.

Finally, condition on all diagonal entries in the original model.
Under this conditional law, the entries above the diagonal are independent with
their original distributions.
Since the bound for fixed diagonal entries is uniform in their
values, averaging over their joint distribution proves the theorem.
\end{proof}

\section*{Acknowledgments}

The authors are grateful to Alexander Litvak for his support
and helpful suggestions.

\section*{Statement on the Use of AI}

The main ideas, results, and proof strategy of this paper were
developed by the authors. During the preparation of the manuscript,
ChatGPT and Codex (OpenAI), including the GPT-5.6 Sol and GPT-6 Pro
models, were used as auxiliary tools to check and refine certain
technical details of the block comparison and thin--thick arguments
after the underlying ideas and proof structure had been supplied by
the authors, and to assist with language editing. All AI-assisted
material included in the manuscript was independently reviewed and
verified by the authors, who take full responsibility for the
mathematical content and final text.

\bibliographystyle{amsplain}
\bibliography{ref}
\end{document}